\documentclass[a4paper,pdftex,reqno,10pt]{amsart}

\allowdisplaybreaks[1]

\usepackage{geometry}

\usepackage{graphicx}
\usepackage{enumerate}
\usepackage{courier}
\usepackage{times}
\usepackage{amsmath,amssymb}

\usepackage{hyperref}

\theoremstyle{plain}
\newtheorem{Theorem}{Theorem}
\newtheorem{Proposition}[Theorem]{Proposition}

\newtheorem{Corollary}[Theorem]{Corollary}
\newtheorem{Lemma}[Theorem]{Lemma}

\newenvironment{Proof}
{\begin{trivlist}\item[]{{\sc Proof.}}}{\hfill{$\square$}\noindent\end{trivlist}}

\theoremstyle{definition}

\theoremstyle{remark}

\DeclareMathOperator{\rk}{rk}

\newcommand{\gaussmset}[2]{\left[\begin{smallmatrix}{#1}\\{#2}\end{smallmatrix}\right]}
\newcommand{\F}{\ensuremath{\mathbb{F}}}
\newcommand{\Fq}{\ensuremath{\F_q}}
\newcommand{\Fqn}{\ensuremath{V}}

\newcommand{\Gqnk}{\ensuremath{\gaussmset{\Fqn}{k}}}

\newcommand{\rdist}{\mathrm{d}_{\mathrm{r}}}
\newcommand{\sdist}{\mathrm{d}_{\mathrm{s}}}

\begin{document}


\title{Subspaces intersecting in at most a point}


 \author{Sascha Kurz}
 \address{Sascha Kurz, University of Bayreuth, 95440 Bayreuth, Germany}
 \email{sascha.kurz@uni-bayreuth.de}

\abstract{We improve on the lower bound of the maximum number of planes in $\operatorname{PG}(8,q)\cong\F_q^{9}$ pairwise 
intersecting in at most a point. In terms of constant dimension codes this leads to $A_q(9,4;3)\ge q^{12}+ 2q^8+2q^7+q^6+2q^5+2q^4-2q^2-2q+1$. 
This result is obtained via a more general construction strategy, which also yields other improvements.
}}

\maketitle

\section{Introduction}
Let $V\cong{\F}_q^v$ be a $v$-dimensional vector space over the finite field {\Fq} with $q$ elements. We call each $k$-dimensional linear subspace 
of $V$ a $k$-space, also using the terms points, lines, and planes for $1$-, $2$-, and $3$-spaces, respectively. Two $k$-spaces $U$, $W$ are said to 
trivially intersect or to be disjoint if $\dim(U\cap W)=0$, i.e., $U$ and $W$ do not share a common point. Sets of $k$-spaces that are pairwise disjoint 
are called partial $k$-spreads, see \cite{ubt_eref42176} for a recent survey on bounds for their maximum possible sizes. In finite projective geometry 
they are a classical topic. Here we study the rather similar objects of sets of $k$-spaces which pairwise intersect in at most a point and have 
large cardinality. More generally, we can use the subspace distance  $\sdist(U,W) = \dim(U+W)-\dim(U \cap W)=\dim(U)+\dim(W)-2\dim(U \cap W)$ 
to define $A_q(v,d;k)$ as the maximum number of $k$-spaces in $\F_q^v$ that have minimum subspace distance $d$, i.e., that intersect in a subspace 
of dimension at most $k-d/2$. Since those sets, which are also called constant dimension codes, have applications in error correcting random network 
coding, see e.g.\ \cite{MR2451015}, bounds for $A_q(v,d;k)$ have been studied intensively in the literature. For the currently best known 
lower and upper bounds we refer to the online tables \url{http://subspacecodes.uni-bayreuth.de} and the associated survey \cite{HKKW2016Tables}. 
Due to this connection, we also call sets of $k$-spaces codes and call their elements codewords. 

Due to combinatorial explosion, it is in general quite hard to obtain improvements for $A_q(v,d;k)$ when the dimension $v$ of the ambient space 
is small, say $v\le 11$. 
Our main motivation for this paper is the recently improved parametric lower bound $A_q(9,4;3)\ge q^{12}+ 2q^8+2q^7+q^6+q^5+q^4+1$, 
see \cite[Theorem 3.13]{cossidente2019subspace}. Here, we give a further improved construction for $A_q(9,4;3)$ and generalize the underlying 
ideas to a more general combination of constant dimension codes. The latter constitutes our main Theorem, see Theorem~\ref{main_thm_new}, which 
allows to conclude also other improved parametric constructions.   

\section{Preliminaries}
For two matrices $U,W\in\mathbb{F}_q^{m\times n}$ we define the rank distance $\rdist(U,W):=\rk(U-W)$. A subset $\mathcal{C}\subseteq \mathbb{F}_q^{m\times n}$ is
called a rank metric code. 

\begin{Theorem}(see \cite{delsarte1978bilinear}) 
  \label{thm_MRD_size}
  Let $m,n\ge d$ be positive integers, $q$ a prime power, and $\mathcal{C}\subseteq \mathbb{F}_q^{m\times n}$ be a rank metric
  code with minimum rank distance $d$. Then, $\# \mathcal{C}\le q^{\max\{n,m\}\cdot (\min\{n,m\}-d+1)}$.
\end{Theorem}
Codes attaining this upper bound are called maximum rank distance (MRD) codes. They exist for all choices of parameters. 
A construction can e.g.\ be described using so-called linearized polynomials, see e.g.~\cite[Section V]{MR2451015}.  
If $m<d$ or $n<d$, then only $\# \mathcal{C}=1$ is possible, which can be achieved by a zero matrix and may be summarized to the single upper bound
$\# \mathcal{C}\le \left\lceil q^{\max\{n,m\}\cdot (\min\{n,m\}-d+1)}\right\rceil$.
Using an $m\times m$ identity matrix as a prefix one obtains the so-called lifted MRD codes.
\begin{Theorem}\cite[Proposition 4]{silva2008rank}
  For positive integers $k,d,v$ with $k\le v$, $d\le 2\min\{k,v-k\}$, and $d$ even, the size of a lifted MRD code in $\Gqnk$ with
  subspace distance $d$ is given by $q^{\max\{k,v-k\}\cdot(\min\{k,v-k\}-d/2+1)}$.
\end{Theorem}

\section{Combining subspaces}
\begin{Theorem}
  \label{main_thm_new}
  Let $\mathcal{C}_1$ be a set of $k$-spaces in $\F_q^{v_1}$ mutually intersecting in at most a point, $\mathcal{C}_1^C$ be a subset of $\mathcal{C}_1$ 
  such that all elements are pairwise intersecting trivially, and $\mathcal{C}_2$ be a set of $k$-spaces in $\F_q^{v_2}$ mutually intersecting in at most 
  a point, where $v_2\ge 2k$ and $\#\mathcal{C}_2\ge 1$. If $\F_q^{v_2}$ admits a $(v_2-k)$-space $S$, such that exactly $\Lambda$ elements 
  of $\mathcal{C}_2$ are contained in $S$ and all others intersect $S$ in at most a point, then 
  $$
    A_q(v_1+v_2-k,2k-2;k)\ge \#\mathcal{C}_1\cdot q^{2(v_2-k)}+\#\mathcal{C}_1^C\cdot \left(\#\mathcal{C}_2-q^{2(v_2-k)}-\Lambda\right)+\Lambda.
  $$  
\end{Theorem}
\begin{Proof}
  We embed $\mathcal{C}_1$ in $\F_q^{v_1+v_2-k}$ and choose a $(v_2-k)$-space $S$ disjoint to the span $\langle\mathcal{C}_1\rangle$. For each $U\in\mathcal{C}_1$ we 
  consider the $v_2$-space $K=\langle U,S\rangle$. If $U\in \mathcal{C}_1^C$, we embed $\mathcal{C}_2$ minus the $\Lambda$ codewords contained in $S$ in $K$ 
  such that the embedding contains the $k$-space $U$ and all codewords intersect $S$ in at most a point. If $U\notin \mathcal{C}_1^C$, we embed a lifted MRD code 
  in $K$ such that the embedding contains the $k$-space $U$ and all codewords are disjoint to $S$. If we additionally add $\Lambda$ codewords inside $S$, then 
  we obtain a set $\mathcal{C}$ of $k$-spaces in $\F_q^{v_1+v_2-k}$ of cardinality $\#\mathcal{C}_1^C\cdot \left(\#\mathcal{C}_3-\Lambda\right)+
  \left(\#\mathcal{C}_1-\#\mathcal{C}_1^C\right)\cdot q^{2(v_2-k)}+\Lambda$, since the matching  lifted MRD code has cardinality $q^{2(v_2-k)}$. For two different 
  $W,W'\in\mathcal{C}$ we have to show that they do intersect in at most a point. By construction, there exist 
  $U,U'\in\mathcal{C}_1$ such that $W\le K:=\langle U,S\rangle$ and $W'\le K':=\langle U',S\rangle$. We have $S\le K\cap K'$ and 
  $v_2-k\le \dim(K\cap K')=v_2-k+\dim(U\cap U')\le v_2-k+1$. If $U=U'$, which we can assume w.l.o.g.\ for $W\le S$ or $W'\le S$, then $\dim(W\cap W')\le 1$. 
  If $U,U'\in \mathcal{C}_1^C$, then $W\cap W'\le S$, so that $\dim(W\cap W')\le 1$. Otherwise we have $\dim(W\cap W'\cap S)=0$, so that also $\dim(W\cap W')\le 1$.      
\end{Proof}

If we choose $v_2=2k$ and $\mathcal{C}_2$ such that there are two disjoint codewords, then $S$ can be chosen as a codeword, i.e., $\Lambda=1$, and all codewords 
except $S$ itself intersect $S$ in at most a point. For brevity, we will calls sets of $k$-spaces that are trivially intersecting and are a subset of a 
some set $\mathcal{C}_1$ of $k$-spaces, a clique.  

\begin{Corollary}
  \label{cor_pavese_cossidente}
  $$A_q(9,4;3)\ge q^{12}+ 2q^8+2q^7+q^6+2q^5+2q^4-2q^2-2q+1$$
\end{Corollary}
\begin{Proof}
  For $k=3$ and $v=6$ we choose $\mathcal{C}_1$ and $\mathcal{C}_2$ as a set of $q^6+2q^2+2q+1$ planes in $\F_q^6$ pairwise intersecting in at most a point 
  \cite[Theorem 2.1]{cossidente2016subspace}. By \cite[Theorem 3.12]{cossidente2019subspace} we can choose a subset $\mathcal{C}_1^C\subseteq \mathcal{C}_1$ of 
  cardinality $q^3-1$.
\end{Proof}
We remark that this improves the very recent lover bound $A_q(9,4;3)\ge q^{12}+ 2q^8+2q^7+q^6+q^5+q^4+1$ \cite[Theorem 3.13]{cossidente2019subspace}. 
As $\mathcal{C}_2$ we might also have chosen the construction from \cite{MR3329980} of the same size.\footnote{The same applies to $\mathcal{C}_1$, i.e., 
we can avoid to use \cite[Theorem 3.12]{cossidente2019subspace}, see the subsequent Footnote~\ref{fn_clique_sizes}.} In our setting we always have $\#\mathcal{C}_1^C\le A_q(6,6;3)=q^3+1$. 
If we replace $\mathcal{C}_2$ in Corollary~\ref{cor_pavese_cossidente} by the set of $q^8+q^5+q^4-q-1$ planes in $\F_q^{7}$ from \cite[Theorem 3]{MR3444245}, then 
the conditions of Theorem~\ref{main_thm_new} are satisfied for $\Lambda=0$ and we obtain
\begin{equation}
  A_q(10,4;3)\ge q^{14}+2q^{10}+2q^9+2q^8+q^7-q^5-2q^4-q^3+q+1.
\end{equation}
However, \cite[Proposition 4.4]{note_linkage} gives a better lower bound.

For a general application of Theorem~\ref{main_thm_new} the presumably hardest part is to analytically determine $\mathcal{C}_1^C$, i.e., a clique in $\mathcal{C}_1$. 
If $\mathcal{C}_1$ itself is obtained via Theorem~\ref{main_thm_new} and a lower bound on the clique size of the corresponding part $\mathcal{C}_2$ is known, then  
can recursively determine suitably large cliques.
\begin{Lemma}
  \label{lemma_clique_induction}
  If $\mathcal{C}$ is obtained from the construction of Theorem~\ref{main_thm_new} and the corresponding part $\mathcal{C}_2$ contains a clique $\mathcal{C}_2^C$ whose 
  elements are disjoint from $S$, then $\mathcal{C}$ admits a subset $\mathcal{C}'$ such that all elements are pairwise intersecting trivially and 
  $\#\mathcal{C}'=\#\mathcal{C}_1^C\cdot \#\mathcal{C}_2^C$.
\end{Lemma}  
\begin{Proof}
  Using the notation from Theorem~\ref{main_thm_new} we construct $\mathcal{C}'$. For each $U\in\mathcal{C}_1^C$ we consider 
  $K:=\langle U,S\rangle$ and choose a clique of cardinality $\#\mathcal{C}_2^C$ in $K$ and add the elements to $\mathcal{C}'$. Using the analysis of the proof of 
  Theorem~\ref{main_thm_new} again and the fact that the elements of $\mathcal{C}'$ all are disjoint to $S$, we conclude that the elements of $\mathcal{C}'$ are pairwise 
  intersecting trivially.  
\end{Proof}

If we choose $\mathcal{C}_2$ according to \cite[Theorem 2.1]{cossidente2016subspace}, we can use \cite[Theorem 3.12]{cossidente2019subspace} to conclude $\#\mathcal{C}_2^C\ge q^3-1$.

\begin{Proposition}
  \label{prop_parametric_series}
  $A_q(6+3t,4;3)\ge \left(q^6+2q^2+2q+1\right)\cdot q^{6t} +\frac{q^{6t}-1}{q^6-1}+\sum\limits_{i=1}^t (2q^2+2q)\cdot \left(q^3-1\right)^i\cdot q^{6(t-i)}$ for all $t\ge 0$.
\end{Proposition}
\begin{Proof}
  For the induction start $t=0$ we choose $\mathcal{C}^{(0)}$ as a set of $q^6+2q^2+2q+1$ planes in $\F_q^6$ pairwise intersecting in at most a point 
  according to \cite[Theorem 2.1]{cossidente2016subspace}, which admits a clique of cardinality $q^3-1$. For the induction step $\mathcal{C}^{(i)}\to\mathcal{C}^{(i+1)}$ 
  we apply Theorem~\ref{main_thm_new} with $v_2=2k$, $\Lambda=1$, $\mathcal{C}_1=\mathcal{C}^{(i)}$, and $\mathcal{C}_2=\mathcal{C}^{(0)}$. By induction, see  
  Lemma~\ref{lemma_clique_induction}, $\mathcal{C}^{(i)}$ admits a clique $\mathcal{C}_1^C$ of cardinality $\left(q^3-1\right)^{i+1}$. The induction 
  hypothesis for the cardinality of $\mathcal{C}^{(i)}$ is
  \begin{equation}
    \label{eq_cardinality_induction}
    \#\mathcal{C}^{(i)}=\left(q^6+2q^2+2q+1\right)\cdot q^{6i} +\frac{q^{6i}-1}{q^6-1}+\sum\limits_{j=1}^i (2q^2+2q)\cdot \left(q^3-1\right)^i\cdot q^{6(i-j)}
  \end{equation}   
  and the induction step, see Theorem~\ref{main_thm_new}, gives $\#\mathcal{C}^{(i+1)}$ as the right hand side of Equation~(\ref{eq_cardinality_induction}), 
  where $i$ is replaced by $i+1$.
\end{Proof}

Another example of a set of planes pairwise intersecting in at most a point, where we can analytically determine a reasonably large clique, is given by 
\cite[Proposition 4.4]{note_linkage}: $A_q(8,4;3)\ge q^{10}+q^6+q^5+2q^4+2q^3+2q^2+q+1$, which is the currently best known lower bound for $q\ge 3$. 
The essential key here is that the code contains a lifted MRD code of cardinality $q^{10}$ for rank distance $2$. By \cite[Lemma 5]{MR3015712} the MRD 
code can be chosen in such a way that it contains a subcode of cardinality $q^5$ and rank distance $3$.\footnote{Using linearized polynomials to described 
the lifted MRD code, a clique of matching size can be described as the set of monomials $ax$ (including the zero polynomial).} Thus we obtain a clique of cardinality $q^5$ and 
can use Theorem~\ref{main_thm_new} with $v_2=6$ and $\Lambda=1$ to conclude
\begin{equation}
  \label{ie_11_4_3_first}
  A_q(11,4;3)\ge q^{16}+q^{12}+q^{11}+2q^{10}+2q^9+2q^8+2q^7+2q^6+1,  
\end{equation}
which strictly improves upon \cite[Proposition 4.4]{note_linkage}. Of course we can iteratively apply the combination with the  $q^6+2q^2+2q+1$ planes in $\F_q^6$ 
to obtain an infinite parametric series as in Proposition~\ref{prop_parametric_series}. The method generalizes to cases where large constant dimension codes 
are obtained by using lifted MRD codes as subcodes, which frequently is the case. Also the constant dimension codes showing 
$A_q(6,4;3)\ge q^6+2q^2+2q+1$ \cite[Lemma 12, Example 4]{MR3329980} and $A_q(7,4;3)\ge q^8+q^5+q^4+q^2-q$ \cite[Theorem 4]{MR3444245} are closely related. They both arise 
by starting from a lifted MRD code, removing some planes, and then extending again with a larger set of planes, cf.~\cite{ai2016expurgation}. Considering just the reduced 
lifted MRD code, we can deduce clique sizes of $q^3-1$ and $q^4$, respectively.\footnote{\label{fn_clique_sizes}Both constructions are stated in the language of linearized 
polynomials. For \cite[Lemma 12, Example 4]{MR3329980} the representation $\F_q^6\cong \F_{q^3}\times \F_{q^3}$ is used and the planes removed from the lifted MRD code 
correspond to $ux^q-u^qx$ for $u\in\F_{q^3}$, so that the monomials $ax$ for $a\in\F_{q^3}\backslash\{\mathbf{0}\}$ correspond to a clique of cardinality $q^3-1$. 
For \cite[Theorem 4]{MR3444245} the representation $\F_q^7\cong W\times \F_{q^4}$, where $W$ denotes the trace-zero subspace of $\F_{q^4}/\F_q$, is used. The planes 
removed from the lifted MRD code correspond to $r\left(ux^q-u^qx\right)$ for $r\in\F_{q^4}\backslash\{\mathbf{0}\}$ and $u\in\F_{q^4}$ with $\operatorname{tr}(u)=1$, so 
that the monomial s $ax$ for $a\in\F_{q^4}$ correspond to a clique of cardinality $q^4$.} If we choose $\mathcal{C}_1$ in Theorem~\ref{main_thm_new} as the mentioned 
code for $A_q(7,4;3)$ and $\mathcal{C}_2$ as the mentioned code for $A_q(6,4;3)$ or the code for $A_q(7,4;3)\ge q^8+q^5+q^4-q-1$, see\cite[Theorem 3]{MR3444245}, then we obtain
\begin{equation}
  A_q(10,4;3)\ge q^{14}+q^{11}+q^{10}+q^8-q^7+2q^6+2q^5+1  
\end{equation}
and    
\begin{equation}
  A_q(11,4;3)\ge q^{16}+q^{13}+q^{12}+q^{10}+q^8-q^5-q^4.
\end{equation}
Both inequalities improve upon the (for $q\ge 4)$ previously best known lower bounds from \cite[Proposition 4.4]{note_linkage} and the latter improves upon 
Inequality~(\ref{ie_11_4_3_first}). 


So, Theorem~\ref{main_thm_new} can yield improved constructions, but of course not all choices of the involved parameters and codes lead to improvements.  
If $v_1<2k$, then $\#\mathcal{C}_1^C\le 1$, so that no strict improvement over known constructions can be 
obtained. For $k>3$ it might be necessary 
to use $v_2>2k$, since no example for $A_q(2k,2k-2;k)>q^{2k}+1$ is known. In \cite{ubt_eref46984} the authors have indeed shown $A_2(8,6;4)=2^{8}+1=257$  
and conjectured $A_q(2k,2k-2;k)=q^{2k}+1$ for all $k\ge 4$. 


In principle it is also possible to generalize Theorem~\ref{main_thm_new} 
to situations where the $k$-spaces can intersect in subspaces of dimension $t$ strictly larger than one. To this end, one may partition $\mathcal{C}_1$ into subsets 
$\mathcal{C}_1^{(0)}$, $\mathcal{C}_1^{(1)}$, \dots, $\mathcal{C}_1^{(t)}$ such that every element from $\mathcal{C}_1^{(i)}$ intersects each different element from 
$\cup_{j=0}^i \mathcal{C}_1^{(j)}$ in dimension at most $i$, which generalizes the partition $\mathcal{C}_1^C$, $\mathcal{C}_1\backslash\mathcal{C}_1^C$. If $S$ is again 
our special subspace and $U\in \mathcal{C}_1^{(i)}$, then codewords in the code in $\langle U,S\rangle$ should intersect $S$ in dimension at most $t-i$, where we may 
also put some additional codewords into $S$. Since we currently have no example at hand that improves upon a best known lower bound for $A_q(v,d;k)$, we refrain 
from giving a rigorous proof and detailed statement.    

\section*{Acknowledgment}
The author would like to thank Thomas Honold for his analysis of possible cliques sizes in the constant dimension codes from \cite[Lemma 12, Example 4]{MR3329980} and  
\cite[Theorem 4]{MR3444245}, see Footnote~\ref{fn_clique_sizes}. The main idea for Theorem~\ref{main_thm_new} is inspired by \cite{cossidente2019subspace}.


\end{document}